\documentclass{scrartcl}

\KOMAoptions{paper=a4}
\KOMAoptions{fontsize = 12pt}
\KOMAoptions{DIV=calc} 

\usepackage[utf8]{inputenc} 
\usepackage[T1]{fontenc} 
\usepackage[ngerman,english]{babel} 
\usepackage{csquotes} 
\usepackage{hyphenat} 
\usepackage{todonotes} 
\usepackage{mathptmx} 
\usepackage{avant}
\usepackage{courier}

\usepackage{enumitem}
\setlist[enumerate]{label*=(\alph*),ref=(\alph*),itemsep=0pt,topsep=5pt}
\setlist[itemize]{itemsep=0pt,topsep=5pt}

\KOMAoptions{DIV=calc}

\usepackage{booktabs} 

\usepackage {amsmath} 
\usepackage{amssymb} 
\usepackage{amscd} 

\newcommand{\CC}{{\mathbf C}}

\newcommand{\NN}{{\mathbf N}}

\newcommand{\RR}{{\mathbf R}}

\newcommand{\ZZ}{{\mathbf Z}}

\newcommand{\TP}{\mathbf{TP}}
\newcommand{\CP}{\mathbf{CP}}
\newcommand{\RP}{\mathbf{RP}}

\newcommand{\FFF}{{\mathcal F}}
\newcommand{\GGG}{{\mathcal G}}
\newcommand{\HHH}{{\mathcal H}}

\newcommand{\PPP}{{\mathcal P}}

\newcommand{\red}{\text{red}}

\DeclareMathOperator{\mult}{mult}
\DeclareMathOperator{\dist}{dist}

\DeclareMathOperator{\Aut}{Aut}

\DeclareMathOperator{\SC}{SC}
\DeclareMathOperator{\SCF}{S}
\DeclareMathOperator{\EI}{E}
\DeclareMathOperator{\RS}{\RR Str}
\DeclareMathOperator{\Stab}{Stab}
\DeclareMathOperator{\Bends}{B}
\DeclareMathOperator{\Fix}{Fix}

\usepackage{graphicx} 
\graphicspath{{pic/}} 

\usepackage[format=plain,labelfont={sf,footnotesize},textfont=small,margin=12pt]{caption} 
\DeclareCaptionLabelSeparator{MySpace}{\enskip } 
\usepackage[style=alphabetic,
						sorting=nyt,
						maxnames=4,
						backend=biber,
						doi=true,
						isbn=false,
						clearlang=true]
						{biblatex}
						
\bibliography{biblio}
\DeclareFieldFormat*{title}{\textit{#1}} 
\renewbibmacro{in:}{} 
\DeclareFieldFormat{journaltitle}{\iffieldequalstr{journaltitle}{ArXiv e-prints}{Preprint}{#1}} 
\DeclareFieldFormat{booktitle}{#1} 
\DeclareRedundantLanguages{english,English}{english,german,ngerman,french}
\DeclareSourcemap{ 
  \maps[datatype=bibtex]{
    \map{
      \step[fieldsource=volume,final]
      \step[fieldset=doi,null]
			\step[fieldset=url,null]
    }  
  }
}

\usepackage{hyperref} 
\usepackage{xcolor} 
\hypersetup{hyperindex,breaklinks} 
\definecolor{darkblue}{RGB}{0,0,170}
\definecolor{darkred}{RGB}{200,0,0}
\hypersetup{citecolor=darkblue, urlcolor=darkblue, linkcolor=darkblue, colorlinks=true} 
\usepackage[all]{hypcap} 

\usepackage{aliascnt} 
\usepackage{amsthm} 

\newtheorem {theorem}{Theorem}[section]

\newaliascnt{proposition}{theorem}
\newtheorem {proposition}[proposition]{Proposition}
\aliascntresetthe{proposition}

\newaliascnt{lemma}{theorem}
\newtheorem {lemma}[lemma]{Lemma}
\aliascntresetthe{lemma}

\newaliascnt{corollary}{theorem}

\aliascntresetthe{corollary}

\newaliascnt{conjecture}{theorem}

\aliascntresetthe{conjecture}

\theoremstyle {definition}

\newaliascnt{definition}{theorem}
\newtheorem {definition}[definition]{Definition}
\aliascntresetthe{definition}

\newaliascnt{example}{theorem}

\aliascntresetthe{example}

\newaliascnt{exercise}{theorem}

\aliascntresetthe{exercise}

\newaliascnt{goal}{theorem}

\aliascntresetthe{goal}

\newaliascnt{construction}{theorem}

\aliascntresetthe{construction}

\theoremstyle {remark}

\newaliascnt{remark}{theorem}
\newtheorem {remark}[remark]{Remark}
\aliascntresetthe{remark}

\newaliascnt{convention}{theorem}

\aliascntresetthe{convention}

\newaliascnt{notation}{theorem}

\aliascntresetthe{notation}

\begin {document}

\title {Lower bounds and asymptotics of real double Hurwitz numbers}
\author {Johannes Rau
{\footnote{{For this work, the author was supported 
by the DFG Research Grant RA 2638/2-1.}}}
{}
}

\maketitle

\begin{abstract}
	\noindent
  We study the real counterpart of double Hurwitz numbers, called real double Hurwitz numbers here.
	We establish a lower bound for these numbers with respect to their dependence on the distribution
	of branch points. 
	We use it to prove, under certain conditions, existence of real Hurwitz covers 
	as well as logarithmic equivalence of real and classical Hurwitz numbers.
	The lower bound is based on the \enquote{tropical} computation of real 
	Hurwitz numbers in \cite{MaRa:TropicalRealHurwitzNumbers}.
\end{abstract}

\section{Introduction}

\label{sec:introduction}

Let $H^\CC_g(\lambda, \mu)$ denote the \emph{complex} (i.e., usual) double Hurwitz numbers.
They count holomorphic maps $\varphi$ from a compact Riemann surface $C$ of genus $g$ to $\CP^1$ 
with $r$ given simple branch points and two additional branch points of ramification profile 
$\lambda, \mu$. Here, $r = l(\lambda) + l(\mu) + 2g - 2$.

A \emph{real structure} $\iota$ for $\varphi : C \to \CP^1$ is an anti-holomorphic involution on $C$ 
such that $\varphi \circ \iota = \text{conj} \circ \varphi$.
The \emph{real} double Hurwitz numbers $H^\RR_g(\lambda, \mu; p)$ count tuples
$(\varphi, \iota)$ of holomorphic maps as above together with a real structure.
Here, we assume that the two special branch points are $0, \infty \in \RP^1$, all branch points lie in $\RP^1$, and
$0 \leq p \leq r$ denotes the number of simple branch points on the positive half axis of $\RP^1 \setminus \{0,\infty\}$.

In this paper, we define numbers $Z_g(\lambda, \mu)$ such that 
\[
	Z_g(\lambda,\mu) \leq H^\RR_g(\lambda,\mu; p) \leq H^\CC_g(\lambda,\mu)
\]
and 
\[
	Z_g(\lambda,\mu) \equiv H^\RR_g(\lambda,\mu; p) \equiv H^\CC_g(\lambda,\mu) \mod 2
\]
for all $0 \leq p \leq r$.
The main results of this paper state that, under certain conditions, 
these lower bounds are non-zero and have
logarithmic asymptotic growth equal to $H^\CC_g(\lambda, \mu)$
(see 
\autoref{prop:sufficientconditionexistence},
\autoref{prop:lowerbound},
\autoref{thm:logasymptotics}).
The definition of $Z_g(\lambda, \mu)$ is based on the tropical computation of real 
double Hurwitz numbers in \cite{MaRa:TropicalRealHurwitzNumbers}. 

Note that the real double Hurwitz numbers $H^\RR_g(\lambda, \mu; p)$ indeed depend on $p$, 
or in other words, on the position of the branch points.
This is the typical behaviour of enumerative problems over $\RR$ (instead of $\CC$). 
It is therefore of interest to find lower bounds for real enumerative problems (with respect to the choice 
of conditions, the branch points here) and use these bounds to prove existence of real solutions
or to compare the number of real and complex solutions of the problem.
Such investigations have been carried out e.g.\ for 
real Schubert calculus \cite{So:RealEnumerative, MuTa}, 
counts of algebraic curves in surfaces passing through points \cite{We, ItKhSh:LogarithmicEquivalence} (see also \cite{GeZi:Construction})
and counts of polynomials/simple rational functions with given critical levels \cite{ItZv, HiRa}. 
In most of these examples, a lower bound is constructed by defining a \emph{signed} count of the real solutions
(i.e., each real solution is counted with $+1$ or $-1$ according to some rule) and showing that this signed count
is invariant under change of the conditions. 
In this paper, we prove similar results for double Hurwitz numbers \emph{without} the explicit
constructions of a signed count.
We hope that this rather simple approach can be extended to other situations
using sufficiently nice combinatorial descriptions of the counting problem.  

One way of defining $Z_g(\lambda, \mu)$ is as follows: It is the number of those 
tropical covers which contribute to the tropical count of $H^\RR_g(\lambda, \mu; p)$ 
with \emph{odd} multiplicity.
We prove in \autoref{thm:zigzaglowerbounds} that these numbers 
provide lower bounds for $H^\RR_g(\lambda,\mu; p)$ as explained above.
Next, we give exact numerical criteria on $\lambda, \mu$ in order for $Z_g(\lambda,\mu)$ to be non-zero,
proving existence of real Hurwitz covers in these case (see \autoref{prop:necessaryconditionexistence}, 
\autoref{prop:sufficientconditionexistence} and \autoref{rem:exactconditionexistence}).  

We  study the asymptotic behaviour of real Hurwitz numbers 
when the degree is increased and only simple ramification points are added.
For example, consider the sequences 
\begin{align*} 
	z_{\lambda, \mu, g}(m) &= Z_g((\lambda, 1^{2m}),(\mu, 1^{2m})), \\
	h^\CC_{\lambda, \mu, g}(m) &= H^\CC_g((\lambda, 1^{2m}),(\mu, 1^{2m})),
\end{align*}
where $(\lambda, 1^{2m})$ stands for adding $2m$ ones to $\lambda$. We prove:

\begin{theorem}[\autoref{thm:logasymptotics}] \label{thm:logasymptoticsIntro}
  Under the existence assumptions for 
	zigzag covers, the two sequences are logarithmically equivalent, 
	\[
		\log z_{\lambda, \mu, g}(m) \sim 4 m \log m \sim \log h^\CC_{\lambda, \mu, g}(m).
	\]
\end{theorem}

This is consistent with the best known results for Welschinger invariants 
and the Hurwitz-type counts of polynomials and rational functions mentioned before.
For better comparison, let us recall the main asymptotic statements from
\cite{ItZv, HiRa}.
Let $S_{\text{pol}}(\lambda_1,\dots,\lambda_k)$ and $S_{\text{rat}}(\lambda_1,\dots,\lambda_k)$ denote the 
signed counts of real polynomials $f(x) \in \RR[x]$ and real simple rational functions $\frac{f(x)}{x-p}$, $f \in \RR[x], p \in \RR$, respectively, 
with prescribed critical levels and ramification profiles as defined in \cite{ItZv, HiRa}.
Set
\begin{align*} 
	s_\text{pol}(m) &= S_{\text{pol}}((\lambda_1, 1^{2m}),\dots,(\lambda_k, 1^{2m})), \\
	s_\text{rat}(m) &= S_{\text{rat}}((\lambda_1, 1^{2m}),\dots,(\lambda_k, 1^{2m})).
\end{align*}
Denote by $h^\CC_\text{pol}(m)$ and $h^\CC_\text{rat}(m)$ the corresponding counts
of complex polynomials/complex rational functions.

\begin{theorem}[{\cite[Theorem 5]{ItZv}}] 
  Assume that each partition $\lambda_i$ satisfies the properties:
	\begin{enumerate}
		\item[(O)] At most one odd number appears an odd number of times in $\lambda_i$. 
		\item[(E)] At most one even number appears an odd number of times in $\lambda_i$. 
	\end{enumerate}
	Then we have 
	\[
	  \log s_\text{pol}(m) \sim  2 m \log m \sim \log h^\CC_\text{pol}(m).
	\]
\end{theorem}

\begin{theorem}[{\cite[Theorem 1.3]{HiRa}}] 
  Assume that each partition $\lambda_i$ satisfies (O) and (E) and
	that $d = |\lambda_i|$ is even (or, an extra parity condition).
	Then 
	\[
	  \log s_\text{rat}(m) \sim  2 m \log m \sim \log h^\CC_\text{rat}(m).
	\]
\end{theorem}

To compare this to our result, note that the non-vanishing assumption for \autoref{thm:logasymptoticsIntro}
is satisfied if both $\lambda$ and $\mu$ satisfy (O). 

\paragraph*{Acknowledgements}

The author would like to thank Renzo Cavalieri, Boulos El Hilany, Ilia Itenberg, Maksim Karev, Lionel Lang and Hannah Markwig
for helpful discussions. 
Parts of this work were carried out at Institut Mittag-Leffler
during my visit of the research program
\enquote{Tropical Geometry, Amoebas and Polytopes}.
Many thanks for the great hospitality and atmosphere!

\section{Real double Hurwitz numbers}

\label{sec:realDHN}

Fix two integers $d > 0$, $g \geq 0$, and two partition $\lambda, \mu$ of $d$.
Set $r := l(\lambda) + l(\mu) + 2g - 2$ and fix a collection of $r$ points 
$\PPP = \{x_1, \dots, x_r\} \subset \CP^1 \setminus \{0,\infty\}$. 

\begin{definition} 
  A \emph{complex ramified cover} of genus $g$, type $(\lambda, \mu)$ and simply branched at $\PPP$
	is a holomorphic maps $\varphi : C \to \CP^1$ of degree $d$ such that
	\begin{itemize}
		\item $C$ is a Riemann surface of genus $g$,
		\item the ramification profiles of $\varphi$ at $0$ and $\infty$ are $\lambda$ and $\mu$, respectively, 
		\item the points in $\PPP$ are simple branch points of $\varphi$.
	\end{itemize}
	Let $\psi : D \to \CP^1$ be another complex ramified cover. 
	An \emph{isomorphism} of complex ramified covers is an isomorphism of Riemann surfaces $\alpha : C \to D$
	such that $\varphi = \psi \circ \alpha$. 
	The \emph{complex double Hurwitz number} 
	\[
	  H^\CC_g(\lambda, \mu) = \sum_{[\varphi]} \frac1{|\Aut^\CC(\varphi)|}
	\]
	is the number of isomorphism classes of complex ramified covers $\varphi$
	of genus $g$, type $(\lambda, \mu)$ and simply branched at $\PPP$,
	counted with weight $1 / \Aut^\CC(\varphi)$.
\end{definition}

\begin{definition} 
	Given a complex ramified cover $\varphi : C \to \CP^1$, a \emph{real structure} on $\varphi$ is a
	antiholomorphic involution $\iota : C \to C$ such that 
	\[
	  \varphi \circ \iota = \text{conj} \circ \varphi.
	\]
	The tuple $(\varphi, \iota)$ is a \emph{real ramified cover}. 
	An isomorphism of real ramified covers 
	$(\varphi : C \to \CP^1, \iota)$ and $(\psi : D \to \CP^1, \kappa)$ 
	is an isomorphism of complex ramified covers $\alpha : C \to D$ such that 
	\[
	  \alpha \circ \iota = \kappa \circ \alpha.
	\]
\end{definition}

\bigskip
As above, fix $d > 0$, $g \geq 0$ and two partition $\lambda, \mu$ of $d$,
and set $r := l(\lambda) + l(\mu) + 2g - 2$. We now fix a collection of $r$ \emph{real} points
$\PPP \subset \RP^1 \setminus \{0,\infty\}$. 
We denote by $\RR_+ \subset \RP^1 \setminus \{0,\infty\}$ the positive 
half of the real projective line and set $p := |\PPP \cap \RR_+|$
the number of positive branch points. 

\begin{definition} 
	The \emph{real double Hurwitz number} 
	\[
	  H^\RR_g(\lambda, \mu; p) = \sum_{[(\varphi, \iota)]} \frac1{|\Aut^\RR(\varphi,\iota)|}
	\]
	is the number of isomorphism classes of real ramified covers $(\varphi, \iota)$
	of genus $g$, type $(\lambda, \mu)$ and simply branched at $\PPP$,
	counted with weight $1/|\Aut^\RR(\varphi,\iota)|$.
\end{definition}

\begin{remark} 
  Let $\varphi : C \to \CP^1$ be a complex ramified cover. Let $\RS(\varphi)$ be the set of real structures for $\varphi$. 
	The group $\Aut^\CC(\varphi)$ acts on $\RS(\varphi)$ by conjugation. The orbits of this action are the isomorphism classes 
	of real ramified covers supported on $\varphi$, and the stabilizer $\Stab(\iota)$ is isomorphic to the group
	of real automorphisms $\Aut^\RR(\varphi,\iota)$. It follows $H^\RR_g(\lambda, \mu; p)$ can be described alternatively as
	\[
	  H^\RR_g(\lambda, \mu; p) = \sum_{[\varphi]} \frac{|\RS(\varphi)|}{|\Aut^\CC(\varphi)|}
	\]
	where the sum runs through all isomorphism classes of \emph{complex} ramified covers. 
	Hence $H^\RR_g(\lambda, \mu; p)$ is the number of complex ramified covers which admit a real structure,
	counted with multiplicity $|\RS(\varphi)| / |\Aut^\CC(\varphi)|$.
	Note also that we can construct an injection $\RS(\varphi) \hookrightarrow \Aut^\CC(\varphi)$
	by fixing any $\iota_0 \in \RS(\varphi)$ and setting $\iota \mapsto \iota \circ \iota_0$. 
	Hence $0 \leq |\RS(\varphi)| / |\Aut^\CC(\varphi)| \leq 1$. 
\end{remark}

\begin{remark} \label{con:noautom}
  Non-trivial automorphisms only occur under rather particular circumstances.
	In fact, it is easy to check that 
	the only complex ramified covers with $\Aut^\CC(\varphi) \neq 1$ are
	$\varphi_d : z \mapsto z^d$ (with $\Aut^\CC(\varphi_d) = \ZZ / d \ZZ$) and 
	compositions $\varphi_d \circ \psi$ where $\psi : C \to \CP^1$ is 
	a hyperelliptic map (here, $\Aut^\CC(\varphi_d \circ \psi) = \ZZ / 2 \ZZ$).
	We can exclude this by assuming $r>0$ and $\{\lambda, \mu \} \not\subset \{(2k), (k, k)\}$.
	It follows $|\RS(\varphi)| = 0,1$, so in this case 
	real and complex double Hurwitz numbers are actual counts of covers
	and, in particular,
	\[
	  0 \leq H^\RR_g(\lambda, \mu; p) \leq H^\CC_g(\lambda, \mu).
	\]
\end{remark}

\begin{remark} 
  For the sake of completeness, let us briefly discuss the cases with non-trivial automorphisms.
	\begin{enumerate}
		\item 
			$g=0, \lambda=\mu=(d)$:
			The only ramified cover in this situation is given by $x \mapsto x^d$. Let $\chi$ be a primitve $d$-th root of unity.
			There are $d$ automorphisms generated by $x \mapsto \chi x$, hence $H^\CC_0((d), (d)) = 1/d$.
			There are also $d$ real structures given by $x \mapsto \chi^i \bar{x}, i=0, \dots, d-1$.
			Hence $H^\RR_0((d), (d); 0) = 1$.
			Note that for $d$ odd all real structures are isomorphic and do not have extra real automorphisms. 
			If $d$ is even, the real structures fall into two isomorphism classes with real automorphism
			$x \mapsto - x$.
		\item 
		  $\{\lambda, \mu \} \subset \{(2k), (k, k)\}$ and $r > 0$:
			There exist exactly $k^{r-1}$ ramified covers with non-trivial automorphisms, and they can 
			be described in the form $C = \{y^2 = f(x)\}, (x,y) \mapsto x^k$. The polynomial $f(x)$ is chosen
			such that, for any simple	branch point $p_i \in \CC \setminus \{0\}$, it has exactly one root 
			among the $k$-th roots of $p_i$. This gives $k^r$ choices for $f$, but $k$ of them are isomorphic
			via $(x,y) \mapsto (\chi x,y)$ where $\chi$ is a $k$-th primitive root. Each of these covers has exactly one extra
			automorphism $(x,y) \mapsto (x,-y)$, hence the total contribution to 
			$H^\CC_g(\lambda, \mu)$ is $k^{r-1} / 2$. \\
			If $k$ is odd, there is exactly one choice to make $f(x)$ real (choosing only real roots),
			and we have two real structures $(x,y) \mapsto (\bar{x},\bar{y})$ and $(x,y) \mapsto (\bar{x},-\bar{y})$
			(alternatively, use $y^2 = \pm f(x)$).
			Hence, this cover contributes $1$ to $H^\RR_g(\lambda, \mu; p)$. 
			If $k$ is even, then $f(x)$ can be chosen real only if $p=0$ or $p=r$. Under this assumption,
			we now have $2^{r-1}$ choices for $f$ (choosing among the two real roots each time, up to switching all of them).
			There are two real structures as before and therefore 
			the total contribution to $H^\RR_g(\lambda, \mu; 0/r)$ is $2^{r-1}$.
	\end{enumerate}
  It follows that $H^\RR_g(\lambda, \mu; p) \leq H^\CC_g(\lambda, \mu)$ holds 
	as long as we exclude:
	\begin{itemize}
		\item $g= 0, \{\lambda, \mu \} \subset \{(d), (\frac{d}{2}, \frac{d}{2})\}$,
		\item $g>0, \{\lambda, \mu \} \subset \{(d), (\frac{d}{2}, \frac{d}{2})\}$, and $d=2$ or $d=4$.
	\end{itemize}
\end{remark}

\section{Tropical double Hurwitz numbers}

\label{sec:tropicalDHN}

In this section we recall the tropical graph counts from 
\cite{CaJoMa:TropicalHurwitzNumbers, BeBrMi:TropicalOpenHurwitzNumbers, GuMaRa, MaRa:TropicalRealHurwitzNumbers}
which compute complex resp.\ real Hurwitz numbers.
We include this summary here for the reader's convenience, using the occasion to 
adapt the definitions to the case of double Hurwitz numbers and introducing
a different convention regarding multiplicities of real tropical covers.

Throughout the following, the word \emph{graph} stands for a finite connected graph $G$ without two-valent vertices. 
The one-valent vertices of $G$ are called \emph{ends}, the higher-valent vertices are \emph{inner vertices}. 
The edges adjacent to an end are called \emph{leaves}, other edges are \emph{inner edges}.
We use the same letter $G$ to denote the the topological space obtained by gluing intervals $[0,1]$ according to the graph structure.
We denote by $G^\circ$ the space obtained by removing the one-valent vertices, called the \emph{inner part} of $G$.
The \emph{genus} of $G$ is the first Betti number $g(G) := b_1(G)$.

A \emph{(smooth, compact) tropical curve} $C$ is a graph together with 
a length $l(e) \in (0, +\infty)$ assigned to any inner edge $e$ of $C$. 
This induces a complete inner metric on $C^\circ$ such that each half-open leaf 
is isometric to $[0, +\infty)$ and each inner edge $e$ is isometric to $[0, l(e)]$. 
There is one exception: The graph $\TP^1$ which consists of a single edge with two one-valent endpoints, 
in which case $C^\circ$ is isometric to $\RR$. 
By use of this construction, the data of lengths $l(e)$ is equivalent to the data of an complete inner metric on $C^\circ$. 

A \emph{piecewise $\ZZ$-linear map} between two tropical curves $C,D$ is a continuous map 
$\varphi : C \to D$ such that for any edge $e \subset C$ and any pair $x,y \in e \cap C^\circ$
there exists $\omega \in \ZZ$ such that 
\[
  \dist(\varphi(x), \varphi(y)) = \omega \dist(x, y).
\]
In particular, $\varphi(C^\circ) \subset D^\circ$. 
By continuity it follows that $\omega=:\omega(e)$ only depends on $e$, we call it the \emph{weight} of $e$ (under $\varphi$).

An \emph{isomorphism} $\varphi : C \to D$ is a bijective piecewise $\ZZ$-linear map whose inverse is also piecewise $\ZZ$-linear.
Equivalently, isomorphisms between $C$ and $D$ can be described by isometries between $C^\circ$ and $D^\circ$.

Given a piecewise $\ZZ$-linear map, by momentarily allowing two-valent vertices, we can find subdivisions
of $C$ and $D$ such that for any edge $e$ of $C$, $\varphi(e)$ is either a vertex or an edge of $D$. 
For any vertex $x$ of $C$ and edge $e'$ of $D$ such that $\varphi(x) \in e'$, we can define
\begin{equation} \label{eq:localdegree} 
  \deg_{e'}(\varphi,x) := \sum_{\substack{e \text{ edge of } C \\ x \in e, \varphi(e) = e'}} \omega(e).
\end{equation}
The map $\varphi$ is a \emph{tropical morphism} if $\text{deg}(\varphi,x) := \text{deg}_{e'}(\varphi,x)$ does not depend
on $e'$, for $x$ fixed. 
In this case, the sum
\[
  \deg(\varphi) := \sum_{\substack{e \text{ edge of } C \\ \varphi(e) = e'}} \omega(e)
\]
is also independent on the choice of $e'$ is called the \emph{degree} of $\varphi$.
Note that isomorphisms are tropical morphisms of degree one.  

Fix two integers $d > 0$, $g \geq 0$, and two partition $\lambda, \mu$ of $d$.
Set $r := l(\lambda) + l(\mu) + 2g - 2$ and fix $r$ points $\PPP \subset \RR$. 
We assume $r > 0$, i.e., we exclude the exceptional case $g=0, \lambda=\mu=(d)$.

\begin{definition} 
  A \emph{tropical cover} of genus $g$, type $(\lambda, \mu)$ and simply branched at $\PPP$ 
	is a tropical morphism $\varphi : C \to \TP^1$ of degree $d$ 
	such that 
	\begin{itemize}
	  \item $C$ is a tropical curve of genus $g$,
		\item $\varphi^{-1}(\RR) = C^\circ$, or equivalently, $f$ is non-constant on leaves,
		\item $\lambda = (\omega(e), e \in L_{-\infty})$ and $\mu =  (\omega(e), e \in L_{+\infty})$, where $L_{\pm \infty}$ denotes the set
					of leaves such that $\pm \infty \in \varphi(e)$, respectively,
		\item each $x \in \PPP$ is the image of an inner vertex of $C$.
	\end{itemize}
	Let $\psi : D \to \TP^1$ be another such tropical cover. An \emph{isomorphism} of the tropical covers
	is an isomorphism $\alpha : C \to D$ of tropical curves such that $\varphi = \psi \circ \alpha$. 
	The \emph{complex multiplicity} of $\varphi$ is 
	\[
	  \mult^\CC(\varphi) := \frac{1}{|\Aut(\varphi)|} \prod_{\substack{e \text { inner} \\ \text{edge of } C}} \omega(e).
	\]
\end{definition}

\bigskip
The following properties are easy to check (see e.g.\ \cite[Section 5]{CaJoMa:TropicalHurwitzNumbers}).
Given $g, (\lambda, \mu)$ and $\PPP$, there is a finite number of isomorphism classes of tropical covers of that type.
Moving the points $x_1, \dots, x_r$ changes the metric structure of these covers, but not their combinatorial 
nor weight structure. 
Moreover, any tropical cover $\varphi : C \to \TP^1$ satisfies the following properties:
\begin{itemize}
	\item The curve $C$ is three-valent, i.e., all vertices are either ends or three-valent.
	\item The map $\varphi$ is non-constant on every edge. In particular, $\mult^\CC(\varphi) \neq 0$.
	\item The curve $C$ has $l(\lambda) + l(\mu)$ ends and $r$ inner vertices, one in the preimage of each $x  \in \PPP$. 
	\item All automorphisms of $\varphi$ are generated by symmetric cycles or symmetric forks of $\varphi$. 
	      A symmetric cycle\fshyp{}fork is a pair of inner edges\fshyp{}leaves, respectively, which share endpoints,
				have equal weights and the same image under $\varphi$. 
\end{itemize}

\begin{theorem}[\cite{CaJoMa:TropicalHurwitzNumbers}] \label{thm:complexcorrespondence}
  The complex Hurwitz number $H^\CC_g(\lambda, \mu)$ is equal to 
	\[
	  H^\CC_g(\lambda, \mu) = \sum_{[\varphi]} \mult^\CC(\varphi),
	\]
	where the sum runs through all isomorphism classes $[\varphi]$ of tropical covers of genus $g$,
	type $(\lambda, \mu)$ and simply branched at $\PPP \subset \RR$. 
\end{theorem}

We will now present the corresponding statement for \emph{real} double Hurwitz numbers. 
It is convenient to use slightly different definitions than in \cite{MaRa:TropicalRealHurwitzNumbers},
see \autoref{rem_differenceOlderPaper} for a comparison.

Let $\varphi : C \to \TP^1$ be a tropical cover.
An edge is \emph{even} or \emph{odd} if its weight is even or odd, respectively.
A symmetric cycle or fork $s$, we denote by $\omega(s)$ the weight of one of its edges, 
and $s$ is \emph{even} or \emph{odd} if $\omega(s)$ is even or odd,
respectively. Moreover, will use the following notation.
\begin{itemize}
	\item $\SCF(\varphi)$ is the set of symmetric cycles and symmetric \emph{odd} forks.
	\item $\SC(\varphi) \subset \SCF(\varphi)$ is the set of symmetric cycles.
	\item For $T \subset \SCF(\varphi)$, $C \setminus T^\circ$ is the subgraph of $C$ obtained
				by removing the interior of the edges contained in the cycles\fshyp{}forks of $T$. 
	\item $\EI(\varphi)$ is the set of even inner edges in $C \setminus \SCF(\varphi)^\circ$,
				i.e., those which are not contained in a symmetric cycle. 
\end{itemize}

	\begin{figure}[]
		\centering
		\input{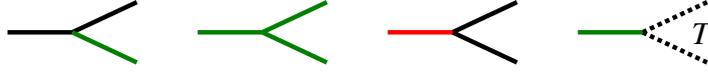}
		\caption{The four types of positive vertices. Odd edges are drawn in black, even edges in colours.
Dotted edges are part of a symmetric fork or cycle contained in $T$.}
		\label{PosVertices}
	\end{figure}%

	\begin{figure}[]
		\centering
		\input{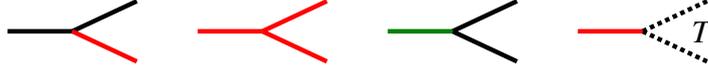}
		\caption{The four types of negative vertices.}
		\label{NegVertices}
	\end{figure}%

\begin{definition} 
	Let $\varphi : C \to \TP^1$ be a tropical cover.
	A \emph{colouring} $\rho$ of $\varphi$ consists of a subset $T_\rho:=T \in \SCF(\varphi)$ and 
	the choice of a colour \emph{red} or \emph{green} for each component of the subgraph of
	even edges of $C \setminus T^\circ$. The tuple $(\varphi, \rho)$ is a \emph{real tropical cover}. 
	An \emph{isomorphism} of real tropical covers is an isomorphism of tropical covers 
	which respects the colouring. 
	The \emph{real multiplicity} of a real tropical cover is 
	\begin{equation} \label{eq:realmult} 
	  \mult^\RR(\varphi, \rho) = 2^{|\EI(\varphi)| - |\SCF(\varphi)|} \prod_{s \in \SC(\varphi)} \mult_T(s),
	\end{equation}
	where 
	\begin{equation} \label{eq:multsymmetriccyclefork} 
		\mult_T(s) := 
			\begin{cases}
				\omega(s) & s \in T, \\
				4 & s \notin T, s \text{ even}, \\
				1 & s \notin T, s \text{ odd}.		
			\end{cases}
	\end{equation}
	Given a real tropical cover, a branch point $x_i \in \PPP$ is \emph{positive}
	or \emph{negative} if it is the image of a three-valent vertex as displayed
	in \autoref{PosVertices} or \autoref{NegVertices}, respectively, up to 
	reflection along a vertical line. This induces a splitting of
	$\PPP = \PPP_+ \sqcup \PPP_-$ into positive and negative branch points.
\end{definition}

\bigskip
Fix $d,g$, $(\lambda, \mu)$ and $\PPP \subset \RR$ as before. 
We now additionally fix a sign for each point in $\PPP$. In other words, 
we fix $0 \leq p \leq r$ and choose a splitting 
$\PPP = \PPP_+ \sqcup \PPP_-$ with $|\PPP_+| = p$.

\begin{theorem}[\cite{MaRa:TropicalRealHurwitzNumbers}] \label{thm:realcorrespondence}
  The real Hurwitz number $H^\RR_g(\lambda, \mu; p)$ is equal to 
	\[
	  H^\RR_g(\lambda, \mu; p) = \sum_{[(\varphi,\rho)]} \mult^\RR(\varphi,\rho),
	\]
	where the sum runs through all isomorphism classes $[(\varphi, \rho)]$ of real tropical covers of genus $g$,
	type $(\lambda, \mu)$, and with positive and negative branch points given by $\PPP = \PPP_+ \sqcup \PPP_-$. 
\end{theorem}

\begin{remark} \label{rem_differenceOlderPaper}
  The present definitions differ from \cite{MaRa:TropicalRealHurwitzNumbers} where 
	$T$ was allowed to contain even symmetric forks as well. Since this choice does not 
	affect the sign of the adjacent vertex (see second and forth vertex in 
  \autoref{PosVertices} resp.\ \autoref{NegVertices}),
	this leads to a factor of $2$ in the number of possible colouring for each such fork. 
	We compensate this by multiplying the real multiplicity from \cite{MaRa:TropicalRealHurwitzNumbers}
	by the same factor. This follows from 
	\[
	  |\Aut(\varphi)| = 2^{|S(\varphi)| + k},
	\]
  where $k$ is the number of even symmetric fork.
	The present convention describes somewhat larger packages 
	of real ramified covers and is more convenient in the discussion that follows.
\end{remark}

\section{Zigzag covers}

\label{sec:zigzagcovers}

In this section we will focus on real tropical covers with odd multiplicity 
$\mult^\RR(\varphi)$ and use them to establish a lower bound
for the numbers $H^\RR_g(\lambda, \mu; p)$ for $0 \leq p \leq r$.
We restict our attention to the automorphism-free case, i.e., we assume
$r>0$ and $\{\lambda, \mu \} \not\subset \{(2k), (k, k)\}$ from now on
(cf.\ \autoref{con:noautom}).

The philosophy behind our approach is as follows:
On one hand, in real enumerative geometry, lower bounds for the counts in question are typically established
by introducing a \emph{signed} count and showing that this alternative count is invariant 
under change of the continuous parameters of the given problem (as long as chosen generically), see e.g.\ \cite{We,ItZv}. 
On the other hand, a real tropical cover $\varphi : C \to \TP^1$ corresponds to a certain package
of $\mult^\RR(\varphi, \rho)$ real ramified covers. Therefore, from the tropical point of view, the easiest
conceivable notion of signs is to ask for maximal cancellation in the tropical packages,
meaning that the signed count of the package of ramified covers tropicalizing to $\varphi$
is 
\[
  \begin{cases}
		0 & \mult^\RR(\varphi, \rho) \text{ even}, \\
		1 & \mult^\RR(\varphi, \rho) \text{ odd},
	\end{cases}
\]
respectively.
In the following, we show that this approach indeed provides an invariant 
and analyse under which conditions this lower bound is non-trivial.
The first step is to prove some properties of $\mult^\RR(\varphi, \rho)$ and express
the condition of having odd multiplicity in combinatorial terms.

\begin{definition} 
  Let $\varphi : C \to \TP^1$ be a tropical cover. 
	The tropical cover $\varphi^\red : C^\red \to \TP^1$ obtained 
	by replacing each symmetric cycle or odd fork (i.e., all $s \in \SCF(\varphi)$) 
	by an edge or leaf of weight $2 \omega(s)$, respectively, 
	is the \emph{reduced} tropical cover of $\varphi$.
\end{definition}

In general, $\varphi^\red$ has smaller branch locus $\PPP' \subset \PPP$.
Note that
$\SCF(\varphi^\red) = \emptyset$ and hence 
$\mult^\RR(\varphi^\red) = 2^{|\EI(\varphi^\red)|}$ for any colouring. 

\begin{lemma} \label{multreducedcover}
  For any tropical cover $\varphi$ we have $|\EI(\varphi^\red)| = |\EI(\varphi)| - |\SCF(\varphi)|$.
\end{lemma}

\begin{proof}
  From our assumptions $r>0$ and $\{\lambda, \mu \} \not\subset \{(2k), (k, k)\}$ it follows that
	$C^\red \neq \TP^1$. 
	Therefore, $C^\red$ contains at least one inner vertex and all its edges
	are isometric either to $[0,l]$ or $[0,\infty]$. 
	
	Let $\varphi^\red =: \varphi_0, \varphi_1, \dots, \varphi_n := \varphi$
	denote the sequence of tropical covers obtained from reinserting the 
	cycles\fshyp{}forks of $C$,	one by one. 
	Since in each step a new even inner edge is created, 
	the	difference $|\EI(\varphi_i)| - |\SCF(\varphi_i)|$ is constant throughout the process.
	Moreover $\SCF(\varphi^\red) = \emptyset$ by construction, which proves the claim. 
\end{proof}

\begin{lemma} \label{multintegeroddeven} 
  For any real tropical cover $(\varphi, \rho)$ the multiplicity $\mult^\RR(\varphi,\rho)$ is an integer
	whose parity is independent of the colouring $\rho$.
\end{lemma}

\begin{proof}
  By \autoref{eq:multsymmetriccyclefork} we have $\mult_T(s) \equiv \omega(s) \mod 2$ and, 
	in particular, the parity of $\mult_T(s)$ does no depend on the colouring.
	Moreover $|\EI(\varphi)| - |\SCF(\varphi)| \geq 0$ by \autoref{multreducedcover}.
	Hence both claims follow from the definition of $\mult^\RR(\varphi)$ in \autoref{eq:realmult}.
\end{proof}

Given a tropical curve $C$, a \emph{string} $S$ in $C$ is 
a connected subgraph such that $S \cap C^\circ$ is a closed submanifold of $C^\circ$.
In other words, $S$ is either a simple loop or a simple path with endpoints
in $C \setminus C^\circ$.
 
\hyphenation{sub-graph}
Let $\varphi$ be a tropical cover. 
Note that any connected component of the subgraph of odd edges is a string 
in $C$. This follows from the definition of tropical morphisms, cf.\ \autoref{eq:localdegree}, 
which implies that at each inner vertex of $C$
the number of odd edges is either $0$ or $2$.

\begin{figure}[]
		\centering
		\input{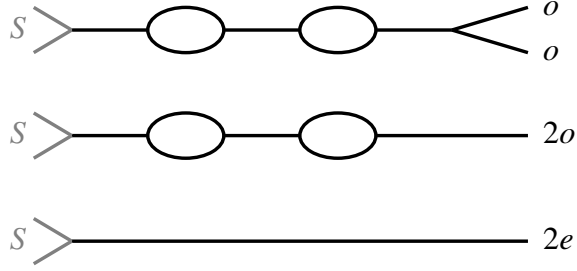}
		\caption{The admissible tails for zigzag covers. The properties of $S$ (turning or not) are not important here.
In the first two cases the number of symmetric cycles can be anything including zero. In the third case, no cycles or forks are allowed.}
		\label{zigzagends}
	\end{figure}

\begin{definition} 
  A \emph{zigzag cover} is a tropical cover $\varphi : C \to \TP^1$ if there exists a 
	subset $S \subset C \setminus \SCF(\varphi)$ such that
	\begin{itemize}
		\item $S$ is either a string of odd edges or consists of a single inner vertex, 
		\item the connected components of $C \setminus S$ are of the form depicted in \autoref{zigzagends}. 
		      Here, all occurring cycles and forks are symmetric and of odd weight. 
	\end{itemize}
\end{definition}

\begin{remark} 
  In \autoref{zigzagends} as well as in the following, we will use the variables $o, o_1, o_2, \dots$ for odd integers
	and $e, e_1, e_2, \dots$ for even integers, respectively.
\end{remark}

\begin{remark}
  Obviously, the set $S$ in the definition of zigzag cover is unique. 
	Moreover, note that the case where $S$ is a single vertex only occurs for special ramification profiles.
	More precisely, we need $\lambda \in \{(e), (o,o)\}$ and 
	$\mu\in \{(e_1,e_2),(e_1, o_2, o_2),(o_1, o_1, o_2, o_2)\}$, or vice versa.
\end{remark}

\begin{proposition} \label{prop:zigzagcoversoddcovers}
  The real tropical cover $(\varphi, \rho)$ is of odd multiplicity if and only if 
	$\varphi$ is a zigzag cover. 
	Moreover, in this case the multiplicity can be expressed as
	\[
	  \mult^\RR(\varphi,\rho) = \prod_{s \in \SC(\varphi) \cap T} \omega(s).
	\]
\end{proposition}

\begin{proof}
  By what has been said so far it follows that $\varphi$ has odd multiplicity if and only if 
	$|\EI(\varphi)| - |\SCF(\varphi)| = 0$ and $\varphi$ does not contain even symmetric cycles.
	Zigzag covers obviously satisfy these properties, so it remains to show the other implication.
	
	By \autoref{multreducedcover} the assumption $|\EI(\varphi)| - |\SCF(\varphi)| = 0$ 
	implies that $\varphi^\red$ does not contain	even inner edges.
	Since $C^\red$ is connected, it follows that there is at most one connected component
	of odd edges in $C^\red$. Therefore $C^\red$ is either an even tripod (three even leaves meeting in one inner vertex) 
	or consists of a single odd component $S$ with some even leaves attached to it. 
	By constrcuction $\varphi$ can be obtained from $\varphi^\red$ by inserting 
	symmetric cycles and symmetric odd forks in the even leaves of $C^\red$.
	Since $C$ does not contain even symmetric cycles, this leads 
	exactly to the three types of tails displayed in \autoref{zigzagends}.  
\end{proof}

\begin{proposition} \label{prop:zigzagworksforallsigns}
  Let $\varphi$ be a zigzag cover simply branched at $\PPP$ and choose
	an arbitrary splitting $\PPP = \PPP_+ \sqcup \PPP_-$ into positive and negative branch points.
	Then there exists a unique colouring $\rho$ of $\varphi$
	such that the real tropical cover $(\varphi, \rho)$ 
	has positive and negative branch points as required.
\end{proposition}

\begin{proof}
	Let $v \in S$ be the vertex from which a given tail $Y$ emanates. 
	Note that if $S = \{v\}$, the vertex cannot be part of a symmetric fork
	since this would imply $\{\lambda, \mu\} = \{(2k), (k,k)\}$.
	Hence, the colour rules from \autoref{PosVertices} and \autoref{NegVertices} impose a unique colouring around $v$ 
	as follows.
	All even edges are coloured in red if $\varphi(v) \in \PPP_-$ and the two edges on the same side of $v$ 
	are both odd (i.e., the string $S$ bends at $v$), or if $\varphi(v) \in \PPP_+$ at least one of the two edges
	on the same side of $v$ is even. In the two opposite cases, all even edges around $v$ are coloured in green.
	
	The next vertex $w$ on $Y$, if it exists, splits the tail into an odd symmetric cycle or fork $s$. 
	We include $s$ in $T$ if the incoming even edge is red and $\varphi(v) \in \PPP_+$ or if 
	the incoming even edge is green and $\varphi(v) \in \PPP_-$. Otherwise, we set $s \notin T$.
	The next vertex $u$ on $Y$, if it exists, closes up an odd symmetric cycle $s$. 
	We colour the outgoing even edge in red if $s \in T$ and $\varphi(v) \in \PPP_+$ or
	of $s \notin T$ and $\varphi(v) \in \PPP_-$. 
	Again, in both cases this process describes the unique local colouring around $v$ compatible 
	with $\PPP = \PPP_+ \sqcup \PPP_-$. 
	Therefore, by continuing this process, 
	we arrive at a compatible colouring of $Y$ and hence all of $C$,
	and uniqueness follows. 
\end{proof}

\begin{definition} \label{def:zigzagnumber}
  The {zigzag number}
	  $Z_g(\lambda,\mu)$
	is the number of zigzag covers of genus $g$, type $(\lambda, \mu)$
	and simply branched at $\PPP \subset \RR$.
\end{definition}

It is easy to check that that $Z_g(\lambda,\mu)$ 
does not depend on the choice of $\PPP \subset \RR$,
see e.g.\ \autoref{rem:combconstructionofcovers}.

\begin{theorem} \label{thm:zigzaglowerbounds}
  Fix $g$, $(\lambda, \mu)$ and $0 \leq p \leq r$ as before. 
	Then the number of real ramified covers is bounded from below by the number 
	of zigzag covers and they have the same parity, i.e.\
	\[
	  Z_g(\lambda,\mu) \leq H^\RR_g(\lambda,\mu; p) \leq H^\CC_g(\lambda,\mu),
	\]
	\[
	  Z_g(\lambda,\mu) \equiv H^\RR_g(\lambda,\mu; p) \equiv H^\CC_g(\lambda,\mu) \mod 2.
	\]
\end{theorem}

\begin{proof}
	The statements involving $Z_g(\lambda,\mu)$ and $H^\RR_g(\lambda,\mu; p)$ 
	follow from \autoref{thm:realcorrespondence} in addition
	with \autoref{prop:zigzagcoversoddcovers} and \autoref{prop:zigzagworksforallsigns}.
	The inequality  $H^\RR_g(\lambda,\mu; p) \leq H^\CC_g(\lambda,\mu)$ is explained
	in  \autoref{con:noautom}, while $Z_g(\lambda,\mu) \equiv H^\CC_g(\lambda,\mu) \mod 2$
	is provided in \autoref{rem:complexmultodd} below for better reference.
\end{proof}

\begin{remark} \label{rem:complexmultodd}
  The first part of \autoref{prop:zigzagcoversoddcovers} holds analogously for 
	complex multiplicities: A tropical cover $\varphi$ is of odd multiplicity $\mult^\CC(\varphi)$ 
	if and only if $\varphi$ is a zigzag cover. It follows that 
	$Z_g(\lambda,\mu) \equiv H^\CC_g(\lambda,\mu) \mod 2$ by 
	\autoref{thm:complexcorrespondence}.
	
	To prove the claim, we copy the proof of \autoref{prop:zigzagcoversoddcovers} 
	replacing (momentarily) $\varphi^\red$ by the	``full reduction'' $\varphi^{\red'}$ 
	in which also \emph{even} symmetric forks	are removed.
	Note that when reinserting a symmetric cycle or fork,
	the complex multiplicity changes by $\omega(s)^3$ or $\omega(s)$, respectively
	(a factor $2$ is cancelled by the automorphism).
	Hence, if $\mult^\CC(\varphi)$ is odd, $C$ does not contain even symmetric cycles and forks
	and $\varphi^{\red'}$ does not contain even inner edges. The remaining argument is as before.
\end{remark}

\begin{remark} 
  This result is \enquote{optimal} in the following sense.
	\begin{itemize}
		\item In principle, we could count the zigzag covers with their multiplicities as calculated in
		      \autoref{prop:zigzagcoversoddcovers}, but these multiplicities do depend on the colouring
					and hence on the signs $\PPP = \PPP_+ \sqcup \PPP_-$.
					In particular, there is one choice of signs for which $T = \emptyset$ and hence
					$\mult(\varphi,\rho) = 1$.
		\item There are no other covers which contribute to any sign distribution, as the next lemma shows.
	\end{itemize}
\end{remark}

\begin{remark} \label{rem:RefinedInvariants}
  It is possible to define (several) \emph{refined} invariants $R_g(\lambda, \mu) \in \ZZ[q^{\pm}]$ 
	in the sense of \cite{BlGoeFockSpaces,BlGoe}. These are counts of tropical covers with polynomial multiplicities
	such that the specializations 
	\begin{align*} 
		R_g(\lambda, \mu)(1) &= H^\CC_g(\lambda,\mu), \\
		R_g(\lambda, \mu)(-1) &= Z_g(\lambda,\mu)
	\end{align*}
	hold. To understand the properties of these refined counts is work in progress together with 
	Boulos El Hilany and Maksim Karev.
\end{remark}

\begin{lemma} 
  Let $\varphi$ be a tropical cover which admits a colouring compatible with 
	$\PPP = \PPP_+ \sqcup \PPP_-$
	for all possible splittings
	$\PPP_+ \sqcup \PPP_-$. Then $\varphi$ is a zigzag cover.
\end{lemma}

\begin{proof}
  Let us consider of yet another version of reducing $C$, denoted by $\varphi^{\red''}$,
	where only odd \emph{odd} symmetric cycles and odd symmetric forks are removed. 
	By our previous considerations, zigzag covers can be equivalently described by the property
	that $\varphi^{\red''}$ does not contain even inner edges. 
	
  Let $\varphi$ be a non-zigzag cover and let $e$ be an 
	even inner edge $e$ of $C^{\red''}$ with endpoints $v_1, v_2$.
	In $C$, the edge $e$ corresponds to a sequence of
	even edges and odd symmetric cycles. 
	We claim that the signs at the branch points of this sequence cannot 
	be chosen independently.
	Indeed, fix the signs for all branch points except, say, $\varphi(v_2)$.
	Then the same process as in the proof of \autoref{prop:zigzagworksforallsigns},
	starting at $v_1$, shows that there is a unique colouring of the sequence 
	compatible with the chosen signs. In particular, the sign $\varphi(v_2)$ is already
	determined by this data. Hence $\varphi$ does not satisfy the condition of the statement,
	which proves the claim.
\end{proof}

	\begin{figure}[]
		\centering
		\input{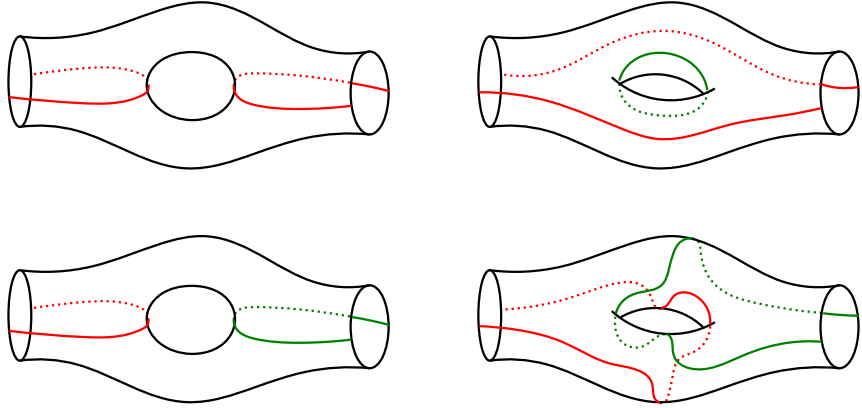}
		\caption{The four possible real structures corresponding of a symmetric cycle (up to switching colours). Here,
the green and red curves represent parts of $\Fix(\iota)$ which are mapped to $\RR_+$ and $\RR_-$, respectively.}
		\label{CycleLift}
	\end{figure}%

\begin{remark} 
  Recall that a Riemann surface with real structure $(C,\iota)$ is called \emph{maximal} if 
	$b_0(\Fix(\iota)) = g + 1$. Following the correspondence of real tropical covers
	and (classical) real ramified covers from \cite{MaRa:TropicalRealHurwitzNumbers}, it is easy to check 
	that all real structures obtained from zigzag covers are maximal. 
	Indeed, note that $S$ accounts for two connected components of the real part
	if it is a loop (mapping to $\RR_+$ and $\RR_-$, repsectively), and for one connected component otherwise.
	Additionally, each symmetric
	cycle produces another connected component as displayed in \autoref{CycleLift}.
	Hence all real ramified covers contributing to the zigzag count are maximal. 
	
	It might also be interesting to specialize to the following type of covers: 
	A real ramified cover $(\psi : C \to \CP^1, \iota)$ 
	is called \emph{of weak Harnack type} if it is maximal and if there exists a single component $H \subset \Fix(\iota)$
	such that $\psi^{-1}(\{0,\infty\}) \subset H$ (i.e., the fibers of $0$ and $\infty$ are totally
	real and lie in a single component of $\Fix(\iota)$ --- we do not impose a condition on the order
	of appearance of the ramification points in $H$, however).
	For zigzag covers, this corresponds to only allowing the upper right lifting shown in 
	\autoref{CycleLift} (in particular, $T = \emptyset$ and $\mult^\RR(\varphi,\rho) = 1$).
	This case occurs for example if $p=0,r$ and if all vertices of $S$ are bends.
	We plan to address these questions and possible connections to 
	refined invariants (see \autoref{rem:RefinedInvariants}) in future work.
\end{remark}

\section{Counting zigzag covers}

\label{sec:countingzigzags}

In this section we discuss existence and asymptotic behaviour of zigzag covers.
We start with a simple observation.

\begin{proposition} \label{prop:necessaryconditionexistence}
  If $Z_g(\lambda,\mu) > 0$, then the number of odd elements which appear an odd number 
	of times in $\lambda$ plus the number of elements which appear an odd number of times
	in $\mu$ is $0$ or $2$.
\end{proposition}

\begin{proof}
  It follows immediately from the three types of tails allowed in a zigzag cover
	that this number is at most $2$. Since it is even, the statement follows. 
\end{proof}

	For existence statements we need slightly stronger assumptions. 
It is useful to introduce some notation for partitions first.
The number of parts in a partition is called the \emph{length} and denoted by 
$l(\lambda)$. The sum of the parts is denoted by $|\lambda|$. 
A partition $\lambda = (\lambda_1, \dots, \lambda_n)$ is called \emph{even} or \emph{odd}
if all parts $\lambda_i$ are even or odd numbers, respectively. 
We denote
\begin{align*} 
  2 \lambda &:= (2 \lambda_1, \dots, 2 \lambda_n), \\
	\lambda^2 &:= (\lambda_1, \lambda_1, \dots, \lambda_n, \lambda_n), \\
	(\lambda, \mu) &:= (\lambda_1, \dots, \lambda_n, \mu_1, \dots, \mu_m),
\end{align*}
where $\mu = (\mu_1, \dots, \mu_m)$ is a second partition. 
Any partition $\lambda$ can be uniquely decomposed into
\[
  \lambda = (2 \lambda_{2e}, 2 \lambda_{2o}, \lambda_{o,o}^2, \lambda_0)
\]
such that
\begin{itemize}
	\item the partition $\lambda_{2e}$ is even,
	\item the partitions $\lambda_{2o}$ and $\lambda_{o,o}$ are odd,
	\item the partition $\lambda_{0}$ is odd and does not have any multiple entries. 
\end{itemize}
We call this the \emph{tail decomposition} of $\lambda$.
Note that $l(\lambda_0) \equiv |\lambda| \mod 2$.
In terms of this notation, the necessary condition of \autoref{prop:necessaryconditionexistence}
can be stated as $l(\lambda_0, \mu_0) = 0,2$. 
Note that $l(\lambda_0, \mu_0)$ is even since
$l(\lambda_0, \mu_0) \equiv |\lambda| + |\mu| = 2d \mod 2$.

\begin{proposition} \label{prop:sufficientconditionexistence} 
  If $l(\lambda_0, \mu_0) \leq 2$ and $(\lambda_{o,o}, \mu_{o,o}) \neq \emptyset$, then 
	there exist zigzag covers of that type, i.e.\ $Z_g(\lambda,\mu) > 0$.
\end{proposition}

\begin{remark} \label{rem:combconstructionofcovers}
  In the following proof, we use a more combinatorial description of tropical covers.
	Let $C$ be a graph of genus $g$ with only one- and three-valent vertices.
	Fix a orientation on $C$ with no oriented loops and pick positive integer weights 
	for the edges of $\Gamma$ such that the balancing condition holds (i.e., 
	for each inner vertex, the sum of outgoing weights is equal to the sum
	of incoming weights). 
	Finally, fix a set $\PPP = \{x_1 < \dots < x_r\} \subset \RR$, where $r$
	is the number of inner vertices of $C$.
	Then for any choice of total order $v_1, \dots, v_r$ on
	the inner vertices of $C$ extending the partial order induced
	by the orientation,
	there exists a unique
	tropical cover $\varphi : C \to \TP^1$ 
	such that for each edge the orientation agrees with the orientation
	of $\TP^1$ (from $-\infty$ to $+\infty$) under $\varphi$, the given weights agree with the ones induced by $\varphi$,
	and  $\varphi(v_i) = x_i$.
\end{remark}

\begin{proof}
  We distinguish the two cases $l(\lambda_0, \mu_0) = 0$ and $l(\lambda_0, \mu_0) = 2$. 
	
	Case $l(\lambda_0, \mu_0) = 2$: We first construct the underlying abstract graph of a zigzag cover. 
	We start from a string graph $S$ and attach tails to $S$: A tail of the first, second or third type
	for each part of $(\lambda_{o,o}, \mu_{o,o})$, $(\lambda_{2o}, \mu_{2o})$ or $(\lambda_{2e}, \mu_{2e})$,
	respectively. The order of appearance of these tails on $S$ can be chosen arbitrarily. Moreover,
	since $(\lambda_{o,o}, \mu_{o,o}) \neq \emptyset$, there is at least one tail of the first type, on which 
	we place $g$ balanced cycles. We obtain a graph $C$ of genus $g$.	
	and the next step is to equip $C$ with an orientation and weights. 
	
	By construction the leaves of $C$ are labelled by parts of $\lambda$ and $\mu$ (the leaves
	of $S$ are assigned to the parts of $(\lambda_0, \mu_0)$). We use each part as weight
	of the corresponding leaf. Moreover, leaves associated to $\lambda$ are oriented towards
	the inner vertex while leaves associated to $\mu$ are oriented towards the end.
	By the balancing condition, there is a unique extension of the orientation and the weight function to all of $C$. 
	Note that indeed all these weights on $S$ turn out to be odd (and, in particular, non-zero).
	Then, by \autoref{rem:combconstructionofcovers}, any choice of total order on the inner vertices
	extending the orientation order gives rise to a zigzag cover $\varphi : C \to \TP^1$.
	
	Case $l(\lambda_0, \mu_0) = 0$: If $l(\lambda_{o,o}, \mu_{o,o}) > 1$, we proceed as before with the following 
	changes. We remove an arbitrary part $\alpha$ from $(\lambda_{o,o}, \mu_{o,o})$ and use it as weight 
	for both leaves of $S$ (instead of a tail of weight $(\alpha, \alpha)$). 
	If $l(\lambda_{o,o}, \mu_{o,o}) = 1$, we need to distinguish two subcases. If $g=0$, we can proceed as in the previous step.
	This could leave us with tails only of the third type, but since no balanced cycles have to be added, this is not a problem.
	
	If $g>0$, we instead do the following. We construct the abstract graph $C$ as before, but now we glue the two ends of $S$
	to a single inner edge $e$, obtaining a graph $C'$ with a (closed) string $S'$. We assign weights to the tail leaves as before.
	Then again by the balancing condition any choice of an odd weight $\omega(e)$ (and orientation) for the gluing edge $e$ fixes an orientation and 
	the weights on all of $C'$. Moreover, there is a (finite) range of choices for $\omega(e)$ such that $S'$ does not 
	turn into an oriented loop and we can proceed as before.
\end{proof}

\begin{remark} \label{rem:exactconditionexistence}  
  The exact conditions for existence of zigzag covers are as follows. We have $Z_g(\lambda,\mu) > 0$ if and only
	if $l(\lambda_0, \mu_0) \leq 2$ and none of the following three cases occurs:
	\begin{itemize}
		\item $g=0$, $(\lambda_{o,o}, \lambda_0, \mu_{o,o}, \mu_0) = \emptyset$ and $l(\lambda_{2e}, \lambda_{2o}, \mu_{2e}, \mu_{2o}) > 3$.
		\item $g=1$, $(\lambda_{2o}, \lambda_{o,o}, \mu_{2o}, \mu_{o,o}) = \emptyset$ and $(\lambda_0, \mu_0) \neq \emptyset$.
		\item $g>1$, $(\lambda_{2o}, \lambda_{o,o}, \mu_{2o}, \mu_{o,o}) = \emptyset$.
	\end{itemize}
\end{remark}

Our next goal is to give lower bounds for $Z_g(\lambda,\mu)$.
Let $\varphi$ be a zigzag cover. We denote by 
\\
\begin{tabular}{rp{.8\textwidth}}
  $a_l, a_r$ & the number of tails of type $o,o$ to the left and right, \\
	$b_l, b_r$ & the number of bends/orientation changes in $S$, with the peak of the bend pointing to the left and right, \\
	$c$        & the number of unbent vertices of $S$, \\
  $g_l, g_r$ & the number of symmetric cycles located on tails to the left and right,
\end{tabular}
\\
respectively. If $S$ is a vertex, we use the values $b_l = b_r = 0$, $c=1$.

\begin{definition} \label{def:unmixed}
  A zigzag cover $\varphi$ is \emph{unmixed} if its simple branch points
	$x_1 < \dots < x_s$, grouped in segments of length $a_l, 2g_l, b_l, c, b_r, 2g_r, a_r$,
	occur as images of
	\begin{itemize}
		\item the symmetric fork vertices of tails of type $o,o$ to the left, 
		\item the vertices of symmetric cycles located on tails to the left,
		\item bends of $S$ with peaks to the left,
		\item unbent vertices of $S$,
		\item the symmetric pattern for bends/tails to the right.
	\end{itemize}
\end{definition}

	\begin{figure}[]
		\centering
		\input{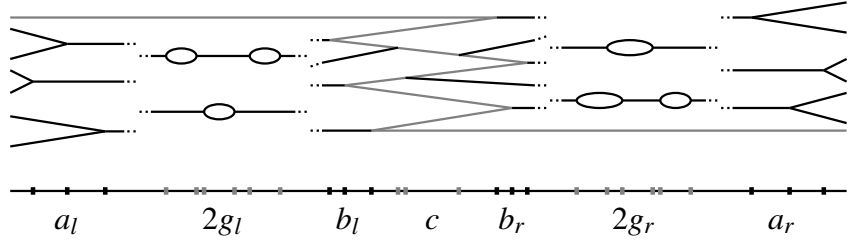}
		\caption{A schematic picture of an unmixed zigzag cover with its various groups of branch points in $\RR$. 
Any permutation of the vertices on top of the branch points belonging to $a_l, b_l, b_r, a_r$, respectively,
gives rise to another unmixed zigzag cover.}
		\label{unmixedzigzag}
	\end{figure}%

\begin{lemma} \label{prop:extendpartialorder}
  Given a weighted and oriented graph $C$ as constructed in \autoref{prop:sufficientconditionexistence},
	there are at least
	\[
	  a_l! \cdot a_r! \cdot b_l! \cdot b_r! 
	\]
	possibilities to turn $C$ into an unmixed zigzag cover $\varphi : C \to \TP^1$.
\end{lemma}

\begin{proof}
	We are interested in the number of total orders on the vertices of $C$ 
	which extend the partial order given by the orientation and,
	when grouped from least to greatest in segments of length $a_l, 2g_l, b_l, c, b_r, 2g_r, a_r$,
	produces the groups of vertices described in \autoref{def:unmixed}.
	It is obvious that such orders exist. Moreover, since
	the partial order restricts to the empty order on the subgroups of vertices
	corresponding to $a_l, b_l, b_r, a_r$, respectively, any permutation
	on these subgroups provides another valid ``unmixed'' order.
	This proves the statement.
\end{proof}

If $\lambda = (\lambda_1, \lambda_2)$, we write $\lambda_1 = \lambda \setminus \lambda_2$.
We write $k \in \lambda$ if $k$ is a part of the partition $\lambda$.
Given a number $k \in \NN$ and two partitions $\lambda, \mu$, consider all (finite) sequences
whose first entry is $k$ and all subsequent entries are obtained by either adding a part 
of $\lambda$ or subtracting a part of $\mu$, using each part exactly once.
We denote by $\Bends(k, \lambda, \mu)$ the maximal number of sign changes that occur in such a sequence.

\begin{theorem} \label{prop:lowerbound} 
  Fix $\lambda, \mu, g$ such that $l(\lambda_0, \mu_0) \leq 2$ and $(\lambda_{o,o}, \mu_{o,o}) \neq \emptyset$.
	We set $\lambda_{\text{tail}} = (\lambda_{2e}, \lambda_{2o}, \lambda_{o,o})$ and 
	$\mu_{\text{tail}} = (\mu_{2e}, \mu_{2o}, \mu_{o,o})$.
	\begin{enumerate}
		\item If $k \in \lambda_0$, then 
			\[
				Z_g(\lambda,\mu) \geq l(\lambda_{o,o})! \cdot l(\mu_{o,o})! \cdot \lfloor \Bends/2 \rfloor! \cdot \lceil \Bends/2 \rceil!,
			\]
			where $\Bends = \Bends(k, 2 \lambda_{\text{tail}}, 2 \mu_{\text{tail}})$. 
		\item If $(\lambda_0, \mu_0) = \emptyset$, $g=0$ and $k \in \lambda_{o,o}$, then 
			\[
				Z_g(\lambda,\mu) \geq (l(\lambda_{o,o})-1)! \cdot l(\mu_{o,o})! \cdot \lfloor \Bends/2 \rfloor! \cdot \lceil \Bends/2 \rceil!,
			\]
			where $\Bends = \Bends(k, 2 (\lambda_{\text{tail}} \setminus (k)), 2 \mu_{\text{tail}})$. 
		\item If $(\lambda_0, \mu_0) = \emptyset$, $g>0$, then 
			\[
				Z_g(\lambda,\mu) \geq l(\lambda_{o,o})! \cdot l(\mu_{o,o})! \cdot (\Bends/2)!^2,
			\]
			where $\Bends = \Bends(1, 2 \lambda_{\text{tail}}, 2 \mu_{\text{tail}})$. 
	\end{enumerate}
\end{theorem}

\begin{proof}
  The three cases are in correspondence with the three types of constructions in \autoref{prop:sufficientconditionexistence}.
	Here, we may assume by symmetry that $\lambda_0 \neq \emptyset$ whenever $(\lambda_0, \mu_0) \neq \emptyset$
	and $\lambda_{o,o} \neq \emptyset$ whenever $(\lambda_{o,o}, \mu_{o,o}) \neq \emptyset$.
	To recall, we summarize the three cases in terms of $S$:
	\begin{itemize}
		\item The leaves of $S$ are weighted by $(\lambda_0, \mu_0)$.
		\item The leaves of $S$ are weighted by a part $k \in \lambda_{o,o}$,
		\item $S$ is a loop.
	\end{itemize}
	
	In each case, $\Bends$ appearing in the statement 
	is the maximal number of bends we can create in $S$	in the corresponding construction. 
	In the first two cases this is straightforward.
	In the loop case, consider a weighted oriented graph $C$ with maximal number of bends on $S$.
	Among the edges of $S$ choose an edge of minimal weight. 
	Following $S$ in the direction of $e$ we subtract $\omega(e) - 1$ from the weights of the edges oriented coherently,
  and add $\omega(e) - 1$ to the edges oriented oppositely. In this way we obtain a new balanced weight function on $C$
	with the same number of bends, but an edge of weight $1$. Hence the maximal number of bends
	is given by $\Bends(1, 2 \lambda_{\text{tail}}, 2 \mu_{\text{tail}})$.

	Back to all three cases, we pick a graph $C$ reaching the maximal number of bends $\Bends$.
	It follows that $b_l = \lfloor \Bends/2 \rfloor$ and $b_r = \lceil \Bends/2 \rceil$,
	respectively. 
	One should note here that $\Bends$ is even if $S$ is a loop, or otherwise, at least one of the ends of $S$ maps 
	to $-\infty$ by our convention from above.
	Moreover, the number of tails of type $o,o$ is
	$a_l = l(\lambda_{o,o})$ and $a_r = l(\mu_{o,o})$, except for the second case, where the ends of $S$ occupy a pair
	of tail weights and hence $a_l = l(\lambda_{o,o})-1$. The statement then follows from \autoref{prop:extendpartialorder}.
\end{proof}

The lower bounds from \autoref{prop:lowerbound} can be used to derive
statements about the asymptotic growth of the numbers $Z_g(\lambda,\mu)$. 

\begin{definition} 
  Given $g \in \NN$ and partitions $\lambda, \mu$ 
	with $|\lambda| = |\mu|$, we set
	\begin{align*} 
	  z_{\lambda, \mu, g}(m) &= Z_g((\lambda, 1^{2m}),(\mu, 1^{2m})), \\
		h^\CC_{\lambda, \mu, g}(m) &= H^\CC_g((\lambda, 1^{2m}),(\mu, 1^{2m})).
	\end{align*}
\end{definition}
 
\begin{proposition} \label{prop:asymptoticszigzag}
  Fix $g \in \NN$, partitions $\lambda, \mu$ with $|\lambda| = |\mu|$ and assume that $l(\lambda_0, \mu_0) \leq 2$.
	Then there exists $m_0 \in \NN$ such that 
	\[
	  z_{\lambda, \mu, g}(m) \geq (m-m_0)!^4
	\]
	for all $m > m_0$.
\end{proposition}

\begin{proof}
  Let $\lambda', \mu'$ be some even partitions of the same integer and $k$ an odd integer.
	We consider the sequence 
	\[
	  \Bends(m) = \Bends(k, (\lambda', 1^{2m}),(\mu', 1^{2m})).
	\]
	We claim that there exists $m_0 \in \NN$ such that $\Bends(m) \geq 2(m-m_0)$ for $m \geq m_0$. 
	Indeed, for sufficiently large $m_0$ we can assume that there exists a maximal sequence 
	(for $B(m_0)$) containing an entry $\pm 1$.
	For $m = m_0 + 1$ we insert a piece of the form $\pm 1 \to \mp 1 \to \pm 1$ at the position of $\pm 1$, and so on,
	showing that $\Bends(m) \geq 2(m-m_0)$. 
	It follows that $\lfloor \Bends(m)/2 \rfloor, \lceil \Bends(m)/2 \rceil \geq m-m_0$ for $m \geq m_0$. 
	
	Note that
	$l((\lambda, 1^{2m})_{o,o}) = l(\lambda_{o,o}) + m$ and $l((\mu, 1^{2m})_{o,o}) = l(\mu_{o,o}) + m$.
  By use of \autoref{prop:lowerbound} and the previous argument 
	we conclude that 
	\begin{align} \nonumber
	  z_{\lambda, \mu, g}(m) &\geq m! \cdot m! \cdot (m - m_0)! \cdot (m - m_0)! \\
		                       &\geq (m-m_0)!^4 \nonumber
	\end{align}
	for all $m \geq m_0$. 
\end{proof}

\begin{theorem} \label{thm:logasymptotics}
  Fix $g \in \NN$, partitions $\lambda, \mu$ with $|\lambda| = |\mu|$ and assume that $l(\lambda_0, \mu_0) \leq 2$.
	Then $z_{\lambda, \mu, g}(m)$ and $h^\CC_{\lambda, \mu, g}(m)$ are logarithmically equivalent. More precisely,
	we have
	\[
	  \log z_{\lambda, \mu, g}(m) \sim 4 m \log m \sim \log h^\CC_{\lambda, \mu, g}(m).
  \]
\end{theorem}

\begin{remark} 
  Let $d$ and $r$ be the degree and number of simple branch points, respectively, 
	of the covers contributing to $z_{\lambda, \mu, g}(m)$.
	Note that the ratio of growth is $4m \sim 2d \sim r$ and therefore
	\[
	  4m \log m \sim  2d \log d \sim r \log r.
	\]
	The variables $d,r$ are more commonly used in the literature, e.g.\ \cite{ItZv, HiRa}.
\end{remark}

\begin{remark} 
	When choosing all simple branch points
	to lie on the positive half axis (i.e., for $p=r$), 
	the logarithmic growth of real double Hurwitz numbers like $H^\RR_0(\lambda,\mu; r)$ 
	can be computed from \cite[Section 5, e.g.\ Theorem 5.7]{GuMaRa}.
\end{remark}

\begin{proof}
  In consideration of \autoref{thm:zigzaglowerbounds}, it suffices to show that $\log z_{\lambda, \mu, g}(m)$
	grows at least as fast and $\log h^\CC_{\lambda, \mu, g}(m)$ grows at most as fast as $4 m \log m$, respectively.
	
	The estimate regarding $\log z_{\lambda, \mu, g}(m)$ follows from \autoref{prop:asymptoticszigzag}
	since $\log ((m-m_0)!) \sim m \log m$.
	
	The estimate for $\log h^\CC_{\lambda, \mu, g}(m)$ (probably classical) can be deduced from the following
	argument. 	
	Let 
	\[
	  H^\CC_g(d) := H^\CC_g((1^d), (1^d))
	\]
	be the complex Hurwitz numbers associated to covers with only simple branch points. 
	The asymptotics of these numbers is computed in \cite[Equation 5]{DuYaZa}
	as 
	\[
	  H^\CC_g(d) \sim C_g \left( \frac{4}{e} \right)^d d^{2d - 5 + \frac{9}{2} g}.
	\]
	Here,  $C_g$ is a constant only depending on $g$. It follows that
	\[
	  \log H^\CC_g(d) \sim 2 d \log d.
	\]
	We finish by showing $H^\CC_g(\lambda', \mu') \leq H^\CC_g(d)$
	for arbitrary partitions $\lambda', \mu'$ of $d$.
	Then
	\[
	  \log h^\CC_{\lambda, \mu, g}(m) \leq \log H^\CC_g(|\lambda| + 2m) \sim 4 m \log m
	\]
	and the claim follows.
	
	To prove $H^\CC_g(\lambda', \mu') \leq H^\CC_g(d)$ we can use monodromy representations.
	Let $\GGG$ and $\HHH$ be the set of monodromy representations corresponding to $H^\CC_g(d)$
	and $H^\CC_g(\lambda', \mu')$. In particular, 
	$|\GGG| = d! \cdot H^\CC_g(d)$ $|\HHH| = d! \cdot H^\CC_g(\lambda', \mu')$.
	Let $\FFF \subset \GGG$ be the subset of representations
	such that 
	\begin{itemize}
		\item the product of the first $d - l(\lambda')$ transpositions gives a permutation of cycle type $\lambda'$,
		\item the product of the last $d - l(\mu')$ transpositions gives a permutation of cycle type $\mu'$.
	\end{itemize}
	We consider the map $\FFF \to \HHH$ given by using as \enquote{special} permutations the products 
	of transpositions as suggested by the definition of $\FFF$. Since any 
	permutation of cycle type $\lambda'$ and $\mu'$ can be factored into a product of 
	$d - l(\lambda')$ and $d - l(\mu')$ transpositions, respectively, the map $\FFF \to \HHH$ is surjective
	and hence $|\HHH| \leq |\GGG|$. This proves the claim.
\end{proof}

The argument can be adapted to prove analogous statements for different types of asymptotics. 
For example, set
\begin{align*} 
	z'_{\lambda, \mu, g}(m)  &= Z_g((\lambda, 2^{m}),(\mu, 1^{2m})), \\
	z''_{\lambda, \mu, g}(m) &= Z_g((\lambda, 2^{m}),(\mu, 2^{m})), 
\end{align*}
and assume in the $z''$ case that $(\lambda_{o,o}, \lambda_0, \mu_{o,o}, \mu_0) \neq \emptyset$ or $g>0$.
A straightforward adaption of \autoref{prop:asymptoticszigzag} shows the following.

\begin{proposition} \label{prop:asymptoticszigzag2}
  Under the above assumptions, there exists $m_0 \in \NN$ such that 
	\begin{align*} 
		z'_{\lambda, \mu, g}(m)  &\geq (m-m_0)!^3, \\
		z''_{\lambda, \mu, g}(m) &\geq (m-m_0)!^2,
	\end{align*}
	for all $m > m_0$.
\end{proposition}

The corresponding series of complex Hurwitz numbers are denoted by
\begin{align*} 
	h'_{\lambda, \mu, g}(m)  &= H^\CC_g((\lambda, 2^{m}),(\mu, 1^{2m})), \\
	h''_{\lambda, \mu, g}(m) &= H^\CC_g((\lambda, 2^{m}),(\mu, 2^{m})).
\end{align*}

\begin{theorem} \label{thm:logasymptotics2}
  Under the above assumptions, we have
	\begin{align*} 
	  \log z'_{\lambda, \mu, g}(m)  &\sim 3 m \log m \sim \log h'_{\lambda, \mu, g}(m), \\
	  \log z''_{\lambda, \mu, g}(m) &\sim 2 m \log m \sim \log h''_{\lambda, \mu, g}(m).
	\end{align*}
\end{theorem}

\begin{remark} 
  The statements can be unified by
	the observation that the logarithmic growth is equal to 
	\[
	  r \log r
	\]
	for all $z,z',z'',h,h',h''$, where $r$ denotes the number of simple branch points.
\end{remark}

\begin{proof}
  We can proceed exactly as for \autoref{thm:logasymptotics} replacing \autoref{prop:asymptoticszigzag}
	by \autoref{prop:asymptoticszigzag2}. It remains to prove that for fixed $k$ the logarithmic growth
	of 
	\begin{align*} 
	  H'_g(m) := H^\CC_g((1^k, 2^m), (1^{k +2m})), \\
	  H''_g(m) := H^\CC_g((1^k, 2^m), (1^k,2^m))
	\end{align*}
	is bounded by $3 m \log m$ and $2 m \log m$, respectively.
	Adapting the previous argument, let $\GGG, \HHH', \HHH''$ be the sets of monodromy representations corresponding to 
	$H^\CC_g(k+2m)$, $H'_g(m)$ and $H''_g(m)$, respectively. Let $\FFF',\FFF'' \subset \GGG$ denote the subsets
	of representations for which the first $m$ transpositions (the first $m$ and the last $m$ transpositions, respectively)
	are pairwise disjoint. We have surjections $\FFF' \to \HHH'$ and $\FFF'' \to \HHH''$ whose fibers have size
	$m!$ and $(m!)^2$, since the factors in product of $m$ pairwise disjoint transpositions can be permuted freely.
	Hence
	\begin{align*} 
	  \log H'_g(m)  &= \log (H^\CC_g(k+2m) / m!) \sim 4 m \log m - m \log m = 3 m \log m, \\
	  \log H''_g(m) &= \log (H^\CC_g(k+2m) / m!) \sim 4 m \log m - 2m \log m = 2 m \log m,
	\end{align*}
	which proves the claim.
\end{proof}

%

\printbibliography

\end {document}